\documentclass{amsart}

\newtheorem{theorem}{Theorem}[section]
\newtheorem{proposition}[theorem]{Proposition}
\newtheorem{corollary}[theorem]{Corollary}

\numberwithin{equation}{section}

\begin{document}

\title{Integral non-hyperbolike surgeries}

\author{Kazuhiro Ichihara}

\address{College of General Education, 
Osaka Sangyo University,
 3--1--1 Nakagaito, Daito, Osaka 574--8530. }

\email{ichihara@las.osaka-sandai.ac.jp}

\subjclass[2000]{Primary 57M50; Secondary 57M25}

\date{February 5, 2005.}

\keywords{exceptional surgery, integral surgery, hyperbolic knot}

\begin{abstract}
It is shown that a hyperbolic knot in the $3$-sphere 
admits at most nine integral surgeries yielding 
$3$-manifolds which are reducible 
or whose fundamental groups are not infinite word-hyperbolic. 
\end{abstract}
\maketitle

%===============================================
\section{Introduction} 

The well-known Hyperbolic Dehn surgery Theorem 
due to Thurston \cite{Th} says that 
each hyperbolic knot admits only finitely many 
Dehn  surgeries yielding non-hyperbolic manifolds. 
A lot of works have been done 
to study how many, when and on which knots 
such exceptional surgeries can occur. 
See \cite{B} for a survey. 

About the number of exceptional surgeries, 
it is conjectured that they are at most TEN, 
and the knot admitting ten is only 
the figure-eight knot in the $3$-sphere $S^3$. 
See \cite[Problem 1.77]{K} for a detail. 
In \cite{HK}, Hodgson and Kerckhoff achieved 
the first universal upper bound SIXTY. 
In \cite{A} and \cite{L}, Agol and Lackenby independently showed that 
there are at most TWELVE 
surgeries yielding non-hyperbolike $3$-manifolds. 
Here a $3$-manifold is called \textit{non-hyperbolike} (in the sense of Agol) 
if it is reducible or does not have infinite 
word-hyperbolic fundamental group. 
Note that hyperbolic implies hyperbolike, 
equivalently non-hyperbolike implies non-hyperbolic. 
Furthermore if the well-known Geometrization Conjecture is true, 
hyperbolike and hyperbolic become equivalent.

The aim of this paper is to present the following result, 
which gives a new upper bound on the number of non-hyperbolike surgeries 
under the assumption that surgery slopes are integral. 
Remark that 
if a meridian-longitude system for $K$ is fixed, 
then surgery slopes for $K$ are parametrized by $\mathbb{Q} \cup \{ 1/0 \}$ 
in a standard way. See \cite{Rol} for example. 

\begin{theorem}\label{thm1}
Let $K$ be a hyperbolic knot 
in a closed orientable $3$-manifold. 
Suppose that 
an arbitrarily chosen meridian-longitude system for $K$ 
is fixed, and by using this, 
surgery slopes for $K$ are parametrized by $\mathbb{Q} \cup \{ 1/0 \}$. 
If two $3$-manifolds obtained by Dehn surgeries on $K$ 
are both non-hyperbolike, and their surgery slopes are both integers, 
then the distance between the surgery slopes is at most eight. 
This implies that there are 
at most NINE integral non-hyperbolike surgeries for $K$. 
\end{theorem}

In particular case of the figure-eight knot  in $S^3$, 
with the standard meridian-longitude system, 
the exceptional surgery slopes 
are $-4, -3, -2, -1, 0, 1, 2, 3, 4$, and $1/0$. 
The slopes $-4$ and $4$ have distance eight, and 
there are nine integral exceptional surgeries. 
Thus our bound is best possible. 

It should be noted that 
Theorem \ref{thm1} assures that 
if some knot admits only integral non-hyperbolike surgeries, 
then the knot have at most TEN non-hyperbolike surgeries, 
that is, nine integral ones and a trivial one. 

For example, it is remarked in \cite{Wu} that 
some class of arborescent knots admit 
only integral surgeries yielding 
reducible, toroidal or Seifert fibered $3$-manifolds. 
Such $3$-manifolds are sometimes called 
non-hyperbolike in the sense of Gordon, 
and are shown to be actually non-hyperbolike in the sense of Agol. 
Precisely we obtain: 

\begin{corollary}
An arborescent knot admits at most TEN 
surgeries yielding 
reducible, toroidal or Seifert fibered $3$-manifolds, 
unless it is a Montesinos knot of type 
$M (x, 1/p, 1/q)$ or its mirror image, 
where $x \in \{ - 1/2n , -1 \pm 1/ 2n , -2 + 1/ 2n \}$, 
and $p$, $q$ and $n$ are positive integers. \qed
\end{corollary}

%----------------------------------
\section{Proof}

We start with preparing definitions and notations. 

As usual, 
we call an embedded circle in a $3$-manifold a \textit{knot}. 
A knot $K$ is called \textit{hyperbolic} 
if its complement $C_K$ admits 
a complete hyperbolic structure of finite volume.

A \textit{slope} is defined to be the isotopy class 
of an unoriented non-trivial simple closed curve on the torus. 
The \textit{distance} between two slopes means 
the minimal geometric intersection number of their representatives.

The complement of an open tubular neighborhood of $K$ 
is called the \textit{exterior} $E_K$ of the knot. 
Given a knot, a new manifold is obtained by 
taking the exterior of the knot and attaching a solid torus back. 
This operation is called a \textit{Dehn surgery} on the knot. 
When one perform a Dehn surgery on a knot, 
a meridian curve of the attached solid torus determines 
a slope on the peripheral torus $\partial E_K$ of the knot. 
This slope is called the \textit{surgery slope} of the Dehn surgery. 
In the following, 
we will use $K(r)$ to denote 
the $3$-manifold obtained 
by Dehn surgery on $K$ with surgery slope $r$.

From now on, we assume that $K$ is a hyperbolic knot in 
a closed orientable $3$-manifold. 
Since $K$ is hyperbolic, 
the complement $C_K$ is regarded as 
a complete hyperbolic $3$-manifold with single cusp. 
The universal cover of $C_K$ is identified with 
the hyperbolic $3$-space $\mathbb{H}^3$. 
Under the covering projection, 
an equivariant set of horospheres 
bounding disjoint horoballs in $\mathbb{H}^3$ 
descends to a torus embedded in $C_K$, 
which we call a \textit{horotorus}. 
As demonstrated in \cite{Th}, on a horotorus $T$, 
a Euclidean structure is obtained by 
restricting the hyperbolic structure of $C_K$. 
By using this structure, the length of a curve on $T$ can be defined. 
Also $T$ is naturally identified with 
the boundary $\partial E_K$ of the exterior $E_K$ of $K$, 
for the image of the horoballs under the covering projection 
is topologically $T$ times half open interval. 
Thus, for a slope $r$ on $\partial E_K$, 
we can define the \textit{length} of $r$ with respect to $T$ 
as the minimal length of the simple closed curves on $T$ 
corresponding to those on $\partial E_K$ with slope $r$.

The following three results will be used in our proof. 
All notations as above will be still used. 

The next proposition was shown by Agol \cite{A} 
and Lackenby \cite{L}, independently.

\begin{proposition}[{\cite[Theorem 6.2]{A}, \cite[Theorem 3.1]{L}}]\label{prop1}
With respect to some horotorus, 
if the length of a slope $r$ on $\partial E_K$ is greater than $6$, 
then the surgered manifold $K(r)$ is irreducible and 
its fundamental group is infinite word-hyperbolic. \qed
\end{proposition}

Now let us choose the maximal horotorus $T$, 
that is, the one bounding the maximal region with no overlapping interior. 
The next proposition holds for such $T$, which was given in \cite{Ad2}. 
See also \cite{Ad1} and \cite{Ad3}.

\begin{proposition}\label{prop2}
With respect to the maximal horotorus, 
every slope on $\partial E_K$ has the length at least $\sqrt[4]{2}$ 
if $K$ is neither the figure-eight knot nor 
the knot $5_2$ in the knot table \cite{Rol}. \qed
\end{proposition}

The next one was obtained in \cite{CM}, 
which is the key to show that the figure eight knot complement 
has the minimal volume among 
orientable $1$-cusped hyperbolic $3$-manifolds.

\begin{proposition}[{\cite[Proposition 5.8]{CM}}]\label{prop3}
For any hyperbolic knot in a closed orientable $3$-manifold, 
the area of the maximal horotorus must be at least 3.35. \qed
\end{proposition}

\begin{proof}[Proof of Theorem $\ref{thm1}$]
We first assume that $K$ is the figure eight knot in $S^3$. 
In this case, exceptional surgeries are completely understood, 
as noted before, the statement holds. 
Next, in the case that $K$ is the knot $5_2$ in $S^3$, 
it is also shown in \cite{BW} that the statement also holds. 

Now, we consider a hyperbolic knot $K$ in general 
neither the figure eight knot nor the knot $5_2$. 
Let $r_1$ and $r_2$ be slopes having distance $\Delta$ 
such that 
they are both integers with respect to 
the fixed meridian-longitude system for $K$ and 
the surgered manifolds $K(r_1)$, $K(r_2)$ are both non-hyperbolike. 

Take the maximal horotorus $T$ in the complement of $K$. 
Let $\mu$ be the closed geodesic on $T$ 
corresponding to the fixed meridian on the boundary $\partial E_K$ 
of the exterior $E_K$ of $K$. 
Let $\gamma_i$ be the closed geodesic on $T$ 
corresponding to a simple closed curve with slope $r_i$ on $\partial E_K$ 
for $i = 1,2$. 
Up to translations, we may assume that 
$\gamma_1$, $\gamma_2$ and $\mu$ have a common intersection point. 

Consider a component $\widetilde{T}$ of the preimage 
of $T$ in the universal cover $\mathbb{H}^3$ of the complement of $K$. 
Since $T$ has a Euclidean structure, 
$\widetilde{T}$ is identified with the Euclidean $2$-plane $\mathbb{E}^2$. 
On $\widetilde{T}$, 
 the preimage of the common intersection point 
of $\gamma_1$, $\gamma_2$ and $\mu$ 
give a lattice. 
By fixing one of the points, say $O$, 
each primitive lattice point corresponds to 
a slope on $T$, and 
the distance between $O$ and a primitive lattice point 
is equal to the length of the corresponding slope. 

We take lattice points $A$ and $B$ 
such that the paths $OA$ and $OB$ are 
lifts of $\gamma_1$ and $\gamma_2$ respectively, and 
the path $AB$ is projected to $\Delta$ times multiple of $\mu$ on $T$. 
Note that the latter condition can be achieved 
by integrality of $r_1$ and $r_2$. 
Also note that the area of the triangle $OAB$ 
is just the half of $\Delta$ times that of $T$. 
The reason is that the parallelogram spanned by $OA$ and $OB$ 
wraps $T$ under the covering projection $\Delta$ times. 

Then we have; 
\begin{itemize}
\item
the length $\overline{OA}$, $\overline{OB}$ 
of the paths $OA$, $OB$ is at most $6$ by Proposition \ref{prop1}, 
\item
the length $\overline{AB}$  of the path $AB$ is at least $\sqrt[4]{2} \Delta$ 
by Proposition \ref{prop2}, 
\item
the area $Area(OAB)$ of the triangle $OAB$ is greater than 
the half of $3.35 \Delta$ by Proposition \ref{prop3}. 
\end{itemize}

Suppose for a contradiction that $\Delta > 8$. 
Let  $\theta$ be the angle between $OA$ and $OB$ so that $ 0 < \theta < \pi$. 
Then by 
$$ 
\overline{OA}^2 + \overline{OB}^2 \leq 6^2 + 6^2 = 72 < 
64 \sqrt{2} < ( \sqrt[4]{2} \Delta )^2 \leq \overline{AB}^2
$$
the angle $\theta$ is greater than $\pi /2 $.

Now, by the Euclidean cosine law, we have
$$ 
\frac{ \overline{AB}^2 - \overline{OA}^2 + \overline{OB}^2 }{ 2 } = 
- \overline{OA} \cdot \overline{OB} \cos \theta \ , 
$$
and so, 
\begin{equation}\label{eq1}
\frac{ \Delta^2 }{ \sqrt{2} }  - 36 < \frac{ ( \sqrt[4]{2} \Delta )^2 - 6^2 - 6^2 }{2}
< 6 \cdot 6 (- \cos \theta) = 36 | \cos \theta| 
\end{equation}
holds.

On the other hand, by the formula of the area of a Euclidean parallelogram, 
we have
$$
2 Area(OAB) = \overline{OA} \cdot \overline{OB} \sin \theta \ ,
$$
and so, 
\begin{equation} \label{eq2}
3.35 \Delta < 6 \cdot 6 \sin \theta = 36 \sin \theta
\end{equation}

Combining Equations \eqref{eq1} and \eqref{eq2}, we have 
$$
\left( \frac{ \Delta^2 }{ \sqrt{2} }  - 36 \right)^2 + ( 3.35 \Delta )^2 
< 36^2 ( \sin^2 \theta + \cos^2 \theta )=  1296. 
$$
From this, we have
\begin{eqnarray*}
\left( \frac{ \Delta^2 }{ \sqrt{2} }  - 36 \right)^2 
+ ( 3.35 \Delta )^2  - 1296 & < & 0 \\
\frac{ \Delta^4 } {2} + ( 11.2225 - 36 \sqrt{2} ) \Delta^2  & < & 0 \\
\Delta^2  & < & 2 ( 36 \sqrt{2} - 11.2225 )  <  80 \\
\Delta & < & \sqrt{80}  \ < \ 8.95 \ .
\end{eqnarray*}

However $\Delta$ must be an integer, and so we would have $\Delta \leq 8$. 
This implies a contradiction to the assumption that $\Delta > 8$. 
\end{proof}

%----------------------------------
\bibliographystyle{amsplain}

\end{document}